\documentclass{article}

\newtheorem{Step}{Step}
\newcommand{\continuum}{2^{\aleph_0}}
\newcommand{\prob}{{\bf P}}
\newcommand{\scrA}{{\cal A}}

\title{The sum of two measurable functions}

\author{Jan Pachl  \\[-6pt] pachl@acm.org}
\date{December 22, 2005\thanks{Transcribed from the author's manuscript dated April 1980.}}

\begin{document}
\maketitle



{\noindent\large\bf Summary}

Following Weizs\"{a}cker~\cite{Weizsaecker}, we use this notation:
For a complete probability space $ ( \Omega, \Sigma, \prob ) $
and a locally convex space $E$,
denote by $ L^0 ( \Omega, \Sigma, \prob, E ) $ the set of all Borel-measurable
functions $ f : \Omega \rightarrow E $ for which the image measure $ f [ \prob ] $
on $E$ is Radon.

In 1976 E. Thomas asked, in a conversation with the author, whether
$ L^0 ( \Omega, \Sigma, \prob, E ) $ is always closed under addition.
The question is motivated by the observation that some of the results in~\cite{Thomas}
can be proved for functions in $ L^0 ( \Omega, \Sigma, \prob, E ) $.

This note presents an example where $ L^0 ( \Omega, \Sigma, \prob, E ) $ is
\underline{not}
closed under addition.
However, Weizs\"{a}cker~\cite{Weizsaecker} showed that this obstacle is not as serious
as would seem.

\vspace{5mm}
{\noindent\large\bf Terminology}

All measures will be probability measures, i.e. positive and with total mass~1.
Say that $ ( X, \scrA, \mu ) $ is a compact Radon measure space
if $X$ is a compact Hausdorff space, $\scrA$ is a sigma-algebra on $X$ containing
all Borel subsets of $X$ and $\mu$ is a complete measure on $\scrA$ such that
\[
\mu B \; = \; \sup \; \{ \;\; \mu K \; | \; K \subseteq B \;\; \mbox{\rm and}
        \;\; K \;\; \mbox{\rm is compact} \; \}
\]
for every $ B \in \scrA $.

When $ ( X , \scrA ) $ and $ ( B , {\cal B} ) $ are two measurable spaces
(sets with sigma-algebras), denote by $ \scrA \otimes {\cal B} $ the product
sigma-algebra on $ X \times Y $;
this is the smallest sigma-algebra on $ X \times Y $ making both projections
$ \pi_1 : X \times Y \rightarrow X $ and $ \pi_2 : X \times Y \rightarrow Y $
measurable.

\vspace{5mm}
{\noindent\large\bf Example}

The example will be constructed in three steps.

\begin{Step}
Construct a compact Radon measure space $ ( X, \scrA, \mu ) $ such that
$ \mu B = 0 $ whenever $ B \in \scrA $ has cardinality less than or equal to
$ \continuum $.
\end{Step}

{\noindent\bf Construction}
Let $I$ be a set of cardinality $ \continuum $,
let $X$ be the compact space $ \{ 0, 1 \}^I $,
and let $ \mu $ be the standard product measure on $X$ (defined on $\scrA$,
the $\mu$-completion of the Borel sigma-algebra in $X$).
That is, $\mu$ is the product of measures on $ \{ 0, 1 \} $ each of which gives
measure $ \frac{1}{2} $ to $ \{ 0 \} $ and $ \frac{1}{2} $ to $ \{ 1 \} $.

For every subset $J$ of $I$, define an automorphism $T_J$ of $ ( X, \scrA, \mu ) $ by
\[
T_J ( \{ x_i \}_{i \in I} ) \; = \; \{ y_i \}_{ i \in I }
\]
where $ y_i = x_i $ for $ i \in J $ and $ y_i = 1 - x_i $ for $ i \in I \backslash J $.

If $ B \in \scrA $ has cardinality $ \leq \continuum $ then there is a set $ J \subseteq I $
such that $ B \cap T_J ( B ) = \emptyset $.
Indeed, choose an injective map $ \alpha : B \times B \rightarrow I $ and define
\[
J \; = \; \{ \; j \in I \; | \; j = \alpha ( \{ x_i \}_{i\in I} , \{ y_i \}_{i\in I} )
       \;\; \mbox{\rm for} \;\; \{ x_i \} , \{ y_i \} \in B \;\; \mbox{\rm and} \;\;
       x_j \neq y_j \; \} \, .
\]

It follows that for each $ B \in \scrA $ of cardinality $ \leq \continuum $ there is
a sequence of $\mu$-automorphisms $ S_1 , S_2 , S_3 , \ldots $ such that
\[
B_k \cap S_k ( B_k ) = \emptyset \; , k = 1, 2, 3, \ldots \, ,
\]
where
\[
B_1 = B \, , \;\; B_{k+1} = B_k \cup S_k ( B_k ) \; .
\]
Hence there are infinitely many disjoint sets of the same measure as B.
Therefore $ \mu B = 0 $.

\begin{Step}
Construct a compact Radon measure space $ ( X, \scrA, \mu ) $
and a measure $\nu$ on the product sigma-algebra
$ \scrA \otimes \scrA $ such that for the ``diagonal''
\\ $ D = \{  (x,x) | x \in X \}$ and the projections
$ \pi_1 : X \times X \rightarrow X $ and $ \pi_2 : X \times X \rightarrow X $
we have
\begin{description}
\item[{\it(i)}] $ \nu G = 1 $ for every $ G \in \scrA \otimes \scrA $
such that $ G \cup D = X \times X $,
\\ and $ \nu H = 1 $ for every $ H \in \scrA \otimes \scrA $
such that $ H \supseteq D $;
\item[{\it(ii)}] $ \pi_1 [ \nu ] = \mu = \pi_2 [ \nu ] $.
\end{description}
\end{Step}

{\noindent\bf Construction}
Take the $ ( X, \scrA, \mu ) $ constructed in Step 1.
Denote by $\beta : X \rightarrow X \times X $
the map defined by $ \beta ( x ) = ( x,x) $.
We have $ \beta^{-1} (G) \in \scrA $ for each $ G \in \scrA \otimes \scrA $;
let $ \nu G = \mu ( \beta^{-1} ( G ) ) $ for $ G \in \scrA \otimes \scrA $ (
that is, $ \nu = \beta [ \mu] $ ).
Since both $ \pi_1 \circ \beta $ and $ \pi_2 \circ \beta $
are the identity map on $X$, (ii) follows.

If $ H \in \scrA \otimes \scrA $ and $ H \supseteq D $ then,
by the definition of $\nu$, we have $\nu H = 1 $.

Take $ G \in \scrA \otimes \scrA $ such that $ G \cup D = X \times X $.
We have
\[
\nu G \; = \; \inf \; \left\{ \; \sum_{n=1}^{\infty} \nu ( B_n \times C_n ) \;\; \left|
       \;\; B_n \, , C_n \in \scrA \;\; \mbox{\rm and } \;
       \bigcup_{n=1}^{\infty} ( B_n \times C_n ) \supseteq G \; \right. \right\}
\]
(see e.g.~\cite{Halmos}, 13.A);
thus it suffices to show that
\[
\sum_{n=1}^{\infty} \nu ( B_n \times C_n ) \geq 1 \;\; \mbox{\rm whenever} \;\;
\bigcup_{n=1}^{\infty} ( B_n \times C_n ) \supseteq G \; , \;\;
B_n , C_n \in \scrA .
\]
Fix such $B_n$, $C_n$ and let $ V = \bigcup_{n=1}^{\infty} ( B_n \times C_n ) $.
Then $ ( X \times X ) \backslash V \subseteq D $.
We show that the cardinality of $ ( X \times X ) \backslash V $ is at most
$ \continuum $:
If $ (x,x) , (y,y) \in ( X \times X ) \backslash V $ and $ x \neq y $
then there is $n$ such that
$ (x,y) \in B_n \times C_n $ and $ (x,x) \not\in B_n \times C_n $;
hence $x$ and $y$ are separated by $C_n$.
It follows that $ (X \times X) \backslash V $ has at most $\continuum $ points.
Consequently,
$ \beta^{-1} ( ( X \times X ) \backslash V ) $ has at most  $\continuum $ points and
\[
\nu ( ( X \times X ) \backslash V ) \;
= \; \mu ( \beta^{-1} ( ( X \times X ) \backslash V ) = 0
\]
by the property of $\mu$.
Thus
\[
\sum_{n=1}^{\infty} \nu ( B_n \times C_n ) \geq \nu V
\; = \; \nu V + \nu ( ( X \times X ) \backslash V ) \; = \; \nu (X \times X) \; = \; 1 \; .
\]
It follows that $ \nu G = 1 $.

\begin{Step}
Construct a complete probability space $ ( \Omega, \Sigma, \prob ) $,
a locally convex space $E$ and two functions $ f,g : \Omega \rightarrow E $ such that
\begin{description}
\item[{\it(a)}] $ f^{-1} (B), g^{-1} (B) \in \Sigma $ for every Borel set $ B \subseteq E $;
\item[{\it(b)}] the image measures $ f [ \prob ] $ and $ g [ \prob ] $ are Radon;
\item[{\it(c)}] the function $ h = f + g $ has the property $ h^{-1} (0) \not\in \Sigma $.
\end{description}
\end{Step}

Thus, in this example, $ L^0 ( \Omega, \Sigma, \prob, E ) $ is not closed under addition. \\[1cm]

{\noindent\bf Construction}
Take the $ ( X, \scrA, \mu ) $, $\nu$, $D$, $\pi_1$ and $\pi_2$ as in Step 2.

Every compact Hausdorff space $Y$ is a topological subspace of a locally convex space
(e.g. let $C(Y)$ be the Banach space of real-valued continuous functions on $Y$;
then $Y$ embeds canonically into the dual of $C(Y)$ endowed with the w* topology).
Fix such an embedding $ e : X \hookrightarrow E $ of $X$ into a suitable
locally convex space $E$.
Let $ ( \Omega, \Sigma, \prob ) $ be the completion of
$ ( X \times X, \scrA \otimes \scrA , \nu ) $, and let
$ f = e \circ \pi_1 $,
$ g = - \, e \circ \pi_2 $.

Now (a) is obvious, and (b) is true because the measures $f[\prob]$ and $g[\prob]$
are continuous images of the Radon measure $\mu$ (by (ii) in Step 2).
Finally, $h^{-1} (0) = D $ and $ D \not\in \Sigma $ in view of (i) in Step 2;
that proves (c).


\end{document}